\documentclass{article}
\usepackage{amsmath}
\font\BBb=msbm10 at 12pt
\newcommand{\Bbb}[1]{\mbox{\BBb #1}}
\newtheorem{theorem}{Theorem}[section]

\newtheorem{cor}[theorem]{Corollary}

\newtheorem{definition}[theorem]{Definition}
\newtheorem{example}[theorem]{Example}

\newcommand{\be}{\begin{equation}}
\newcommand{\ee}{\end{equation}}
\newcommand{\ben}{\begin{enumerate}}
\newcommand{\een}{\end{enumerate}}

\def\beq{\begin{equation}}
\def\eeq{\end{equation}}

\def\pxi{\frac{\partial }{\partial x^i}}
\def\pyi{\frac{\partial}{\partial y^i}}

\newcommand{\qed}{\hspace*{\fill}Q.E.D.}  

\title{{\bf\LARGE On the Projective Ricci Curvature}}
\author{Zhongmin Shen\footnote{supported in part by NSF China (NSFC No. 11671352) }\;  and Liling Sun\footnote{corresponding author}}
\date{February 24, 2020}

\begin{document}

\maketitle
\begin{abstract}

The notion of the Ricci curvature is defined for sprays on a manifold.  With a volume form on a manifold, every spray can be deformed to a projective  spray. The Ricci curvature of a projective  spray is called the projective Ricci curvature.
In this paper, we introduce the notion of projectively Ricci-flat sprays. We establish a global rigidity result for projectively Ricci-flat sprays with nonnegative Ricci curvature. 
Then we study and characterize  projectively Ricci-flat Randers metrics.

{\bf Keywords:} spray,  Finsler metric, Randers metric,  projective Ricci curvature.

{\bf Mathematics Subject Classification 2010:}  53B40, 53C60

\end{abstract}
\section{Introduction}

In Finsler geometry, there are many important
Riemannian quantities such as the Riemann curvature and the Ricci curvature, etc. and non-Riemannian quantities such as the Berwald curvature and the S-curvature, etc.  The Ricci curvature is defined as the trace of the Riemann curvature. Together with the S-curvature, the Ricci curvature plays an important role in Finsler geometry.  The volume of geodesic balls can be  controlled by the lower bounds of the Ricci curvature and the S-curvature (\cite{S2}\cite{S3}). The volume of geodesic balls can be also controlled by a single bound of  the N-Ricci curvature which is the combination of the Ricci curvature and the S-curvature  (\cite{OS1}).

 Let  ${\bf G}$  be a spray on an $n$-dimensional manifold $M$. Given a volume form $dV$ on $M$,   we can construct a new spray  by
 \[ \hat{\mathbf{G}}:=\mathbf{G}+\frac{2\mathbf{S}}{n+1}\mathbf{Y}.\]
The spray $\hat{\bf G}$ is called the projective spray of $({\bf G}, dV)$.
The projective  Ricci curvature of $({\bf G}, dV)$   is defined as the Ricci curvature of $\hat{\bf G}$,  namely,
\begin{equation}\label{7261}
 {\bf PRic}_{(G, dV)}: = {\bf Ric}_{\hat{G}}.
\end{equation}
By a simple computation, we have the following formula for the projective Ricci curvature:
\begin{equation} {\bf PRic}_{(G, dV)}= {\bf Ric} + (n-1) \Big \{  \frac{{\bf S}_{|0}}{n+1}  +\Big [ \frac{\bf S}{n+1} \Big ]^2 \Big \}.\label{630}
\end{equation}
where ${\bf Ric}={\bf Ric}_{G}$ is the Ricci curvature of  the spray ${\bf G}$, ${\bf S} = {\bf S}_{(G, dV)}$ is the S-curvature of $({\bf G}, dV)$ and ${\bf S}_{|0}$ is the covariant derivative of ${\bf S}$ along a geodesic of ${\bf G}$.
It is known  that $\hat{\mathbf{G}}$ remains unchanged under a projective change of ${\bf G}$  with $dV$ fixed,  thus
${\bf PRic}_{({\bf G}, dV)}={\bf Ric}_{\hat{\bf G}}$ is  a projective invariant of $({\bf G}, dV)$ (\cite{S1}). We make the following

\begin{definition}\label{defPRic}A spray  ${\bf G}$ on an $n$-dimensional manifold $M$ is said to be projectively Ricci-flat if there is a volume form $dV$ on $M$ such that
\[ {\bf PRic}_{(G, dV)}=0.\]
 A Finsler metric $F$ on $M$ is said to be {\it projectively Ricci-flat} if the induced spray ${\bf G}={\bf G}_F$ is projectively Ricci-flat.
\end{definition}

It is easy to see that every projectively flat spray on  $R^n$ is projectively Ricci-flat.  For  projectively equivalent sprays, if one of them is projectively Ricci-flat, then so is the other.

We prove the following

\begin{theorem}
Let  ${\bf G}$ be a spray  on an $n$-dimensional manifold $M$.    ${\bf G}$ is projectively Ricci-flat if and only if   there are  a  volume form $dV$ and   a scalar function $f$ on $M$ such that
\be
{\bf Ric}_{G} =-(n-1)\Big \{\Xi_{|0} + \Xi^2 \Big \},\label{PRicf***}
\ee
where
 $``|"$ is the horizontal covariant derivative with respect to $\mathbf{G}$, $f_{0}=f_{x^{m}}y^{m}$,  $\Xi := \frac{{\bf S}}{n+1}-f_0 $ and ${\bf S}= {\bf S}_{(G, dV)}$.
\end{theorem}

If the spray satisfies that ${\bf S}_{(G, dV)} = (n+1) \phi_0$ for some scalar function $\phi $ on $M$, then     $\Xi =0$  for $f =\phi$. We obtain the following

\begin{cor} \label{cor1.3}
Let ${\bf G}$ be a spray on a manifold $M$.   Suppose that ${\bf Ric}_{G} =0$ and ${\bf S}_{(G, dV)} =(n+1) \phi_0$ for some volume form $dV$ and scalar function $\phi$ on $M$, then
${\bf G}$ is projectively Ricci-flat.
\end{cor}

The condition on the S-curvature, ${\bf S}_{(G, dV)} = (n+1) \phi_0$ being an exact $1$-form,  is actually a condition on the spray, independent of the choice of a particular volume form $dV$. Sprays with this property have vanishing $\chi$-curvature $\chi=0$ (\cite{S4}, Proposition 4.1).  
 We have the following rigidity theorem:

\begin{theorem}\label{thmrigidit}
Let ${\bf G}$ be  a complete spray  on an $n$-manifold $M$. Suppose that ${\bf G}$  is projectively Ricci-flat.  If the Ricci curvature ${\bf Ric}_G \geq 0$, then for any volume form $dV$ on $M$, the S-curvature
${\bf S} = (n+1)\phi_0$ for some scalar function $\phi$ on $M$. In this case, ${\bf Ric}_G =0$.
\end{theorem}

To have a better understanding on projectively Ricci-flat Finsler metrics, we consider  a Randers metric  $F=\alpha+\beta$ on an $n$-dimensional manifold $M$, where $\alpha=\sqrt{a_{ij}y^iy^j}$ is a Riemannian metric and $\beta=b_i y^i$ is a $1$-form with $\|\beta\|_{\alpha} < 1$. Put
\[  s_{ij} := \frac{1}{2} ( b_{i;j} - b_{j;i}),\]
where $``;"$ denotes the covariant derivative with respect to Levi-Civita connection of $\alpha$.  Clearly, $\beta$ is closed if and only if $s_{ij}=0$. 
We prove the following
\begin{theorem}\label{6183}
Let $F=\alpha+\beta$ be a Randers metric on an $n$-dimensional manifold $M$. $F$ is projectively Ricci-flat  if and only if there is a scalar function $h$ on $M$ such that
\begin{eqnarray*}
&&\mathbf{Ric}_{\alpha}=2s_{0m}s^m_{\ 0}+\alpha^{2}s^i_{\ j} s^j_{\ i}-(n-1)[h_{0;0}+(h_0)^{2}]\\
&&s^{m}_{~~0;m}=(n-1) h_{x^m}s^{m}_{~~0},
\end{eqnarray*}
where $\mathbf{Ric}_{\alpha}$ denotes the Ricci curvature of $\alpha$.
\end{theorem}

In general, projectively Ricci-flat Randers metrics are not Ricci-flat, and the S-curvature is not almost isotropic. 

In \cite{CSM}, Cheng-Shen-Ma study the projective Ricci curvature  ${\bf PRic}$. They derive a formula for  the projective Ricci curvature of a Randers metric with respect to the Busemann-Hausdorff volume form  $dV_{BH}$. By this formula, they characterize  Randers metrics with ${\bf RRic}=0$ with respect to $dV_{BH}$.  We should point out that the projective Ricci-flatness of Randers metrics defined in \cite{CSM} is slightly different from ours. Thus the statement of Theorem 1.1 in \cite{CSM} is slightly different from  Theorem \ref{6183}.

\section{Preliminaries}
Let $M$ be an $n$-dimensional manifold and $TM$ the corresponding tangent bundle. We denote by $TM_{0}=TM\setminus\{0\}$ the slit tangent bundle. Local coordinates on the base manifold $M$ will be denoted by $(x^{i})$, while the induced local coordinates on $TM$ or $TM_{0}$ will be denoted by $(x^{i} ,y^{i})$. We call $(x^{i},y^{i})$ the standard local coordinate system in $TM$.

A spray $\mathbf{G}$ on $M$ is a smooth vector field on $TM_{0}$ expressed in a standard local coordinate system $(x^{i},y^{i})$ in $TM$ as follows
\begin{equation}
\mathbf{G}=y^{i}\frac{\partial}{\partial x^{i}}-2G^{i}\frac{\partial}{\partial y^{i}},
\end{equation}
where $G^i=G^{i}(x, y)$ are the local functions on $TM$ satisfying
$
G^{i}(x, \lambda y)=\lambda^{2}G^{i}(x, y) ~\text{for}~\forall\lambda>0.
$
 Put
\begin{equation*}
N^{i}_{j}=\frac{\partial G^{i}}{\partial y^{j}}.
\end{equation*}
These are called the {\it nonlinear connection coefficients} of $\mathbf{G}$.   Set $\omega^{i}=dx^{i} $ and $\omega^{n+i}:=dy^{i}+N^{i}_{j}dx^{j}$.
The connection 1-forms of Berwald connection are given by
\begin{equation*}
\omega^{i}_{~j}:=\Gamma^{i}_{jk}dx^{k},
\end{equation*}
where
\begin{equation*}
\Gamma^{i}_{jk}=\frac{\partial N^{i}_{j}}{\partial y^{k}}=\frac{\partial^{2}G^{i}}{\partial y^{j}\partial y^{k}}.
\end{equation*}
We have
\begin{equation*}
d\omega^{i}=\omega^{j}\wedge\omega^{i}_{~j}.
\end{equation*}
The curvature 2-forms of the Berwald connection are defined by
\begin{equation}
\Omega^{i}_{j}=d\omega^{i}_{j}-\omega^{l}_{j}\wedge\omega^{i}_{l}.
\end{equation}
Express
\begin{equation}\label{d2}
\Omega^{i}_{j}=\frac{1}{2}R^{~i}_{j~kl}\omega^{k}\wedge \omega^{l}-B^{~i}_{j~kl}\omega\wedge \omega^{n+l}.
\end{equation}
The two curvature tensors $R^{~i}_{j~kl}$ and $B^{~i}_{j~kl}$ are called \emph{Riemann curvature tensor} and \emph{Berwald curvature tensor}, respectively.

The Riemann curvature $\mathbf{R}_{y}=R^{i}_{~k}dx^{k}\otimes \frac{\partial}{\partial x^{i}}|_{x}:T_{x}M \rightarrow T_{x}M$ is a family of
linear maps on tangent spaces which is defined by
\begin{equation*}
R^{i}_{~k}=2\frac{\partial G^{i}}{\partial x^{k}}-y^{j}\frac{\partial^{2} G^{i}}{\partial x^{j}\partial y^{k}}+2G^{j}\frac{\partial^{2} G^{i}}{\partial y^{j}\partial y^{k}}-\frac{\partial G^{i}}{\partial y^{s}}\frac{\partial G^{s}}{\partial y^{k}}.
\end{equation*}
Without much difficulty, one can show that
\begin{equation}\label{598}
R^{~i}_{j~kl}=\frac{1}{3}(R^{i}_{~k\cdot l}-R^{i}_{~l\cdot k})_{\cdot j}
\end{equation}
and
\begin{equation}\label{611}
R^{i}_{~k}=y^{j}R^{~i}_{j~kl}y^{l}.
\end{equation}
Here and hereafter, notation $``\cdot"$ denotes the vertical derivatives with respect to $y$.
For instance, $f_{\cdot k}=\frac{\partial f}{\partial y^{k}},~~f_{\cdot k\cdot l}=\frac{\partial^{2}f}{\partial y^{k}\partial y^{l}}$, etc.
By \eqref{598} and \eqref{611}, we see that the two curvature tensors $R^{i}_{~k}$ and $R^{~i}_{j~kl}$ can represent each other.
For this reason, they are all called \emph{Riemann curvature tensor} if there is no confusion.

The well-known Ricci curvature is defined by
\begin{equation}\label{692}
{\bf Ric}:=R^{m}_{~m}=y^{j}R^{~m}_{j~ml}y^{l}.
\end{equation}

Every Finsler metric $F$ on a manifold induces a spray $\mathbf{G}_F=y^{i}\frac{\partial}{\partial x^{i}}-2G^{i}(x,y)\frac{\partial}{\partial y^{i}}$ with
the \emph{geodesic coefficients}
\begin{equation}\label{631}
G^{i}=\frac{1}{4}g^{il}\{[F^{2}]_{x^{k}y^{l}}y^{k}-[F^{2}]_{x^{l}}\},
\end{equation}
where $g^{ij}=(g_{ij})^{-1}$. Therefore, every Finsler space is a special spray space.

 Let $\mathbf{G}$ be a spray on an  $n$-manifold $M$ and $dV=\sigma dx^{1}\cdots dx^{n}$ a arbitrary volume form. The $S$-curvature ${\bf S}= {\bf S}_{(G, dV)} $
is defined by
\begin{equation}\label{621}
\mathbf{S}:=\frac{\partial G^{m}}{\partial y^{m}}-y^{m}\frac{\partial }{\partial x^{m}}(\ln \sigma).
\end{equation}

Using the S-curvature  ${\bf S}={\bf S}_{(G, dV)}$, one can modify the spray ${\bf G}$ to
\[  \hat{\bf G}= {\bf G}+\frac{2{\bf S}}{n+1} {\bf Y},\]
where ${\bf Y}= y^i\pyi$ is the canonical vertical vector field. In local coordinates
\[ \hat{G}^i = G^i -\frac{\bf S}{n+1} y^i.\]

\section{Projective Ricci curvature}

 Let $\mathbf{G}$ be a spray on an $n$-dimensional manifold $M$.
 Let $dV$ and $d\tilde{V}$ be volume forms with $dV = e^{-(n+1)f} d\tilde{V}$, where $f=f(x)$ is a scalar function on $M$.
The S-curvatures ${\bf S}= {\bf S}_{(G, dV)}$ and $\tilde{\bf S}= {\bf S}_{(G, d\tilde{V})}$ are related by
\begin{equation}\label{516}
\mathbf{S}
=\tilde{\mathbf{S}}+(n+1) f_0,
\end{equation}
where $f_0 := f_{x^m}(x)y^m$. 
By (\ref{516}), we see that ${\bf S}_{(G, dV)} =(n+1) \phi_0$ if and only if ${\bf S}_{(G, d\tilde{V})} =(n+1) \tilde{\phi}_0$
with $\phi_0 = \tilde{\phi}_0 + f_0$.

By (\ref{630}),  we have
\begin{eqnarray}
{\bf PRic}_{(G, dV)}& = & {\bf Ric}_G + (n-1) \Big \{  \Big [ \frac{\bf S}{n+1}\Big ]^2 + \frac{{\bf S}_{|0}}{n+1} \Big \},\label{eqP1}\\
{\bf PRic}_{(G, d\tilde{V})} & = &{\bf Ric}_G + (n-1) \Big \{  \Big [ \frac{\tilde{\bf S}}{n+1}\Big ]^2 + \frac{\tilde{\bf S}_{|0}}{n+1} \Big \}, \label{eqP2}
\end{eqnarray}
where $``|"$ is the horizontal covariant derivative with respect to $\mathbf{G}$.
It follows from (\ref{eqP1}) and (\ref{eqP2}) that
\begin{equation}\label{5297}
\mathbf{PRic}_{({\bf G}, d\tilde{V})}=\mathbf{{PRic}}_{({\bf G}, dV)}-(n-1)\Big \{ f_{0|0}-f_0^{2}+\frac{2}{n+1}f_0 {\bf S}\Big \},
\end{equation}

\bigskip

\begin{theorem}\label{5310}
Let  ${\bf G}$ be a spray  on an $n$-dimensional manifold $M$.  The following are equivalent:
\ben
\item[(a)]  ${\bf G}$ is projectively Ricci-flat,
\item[(b)]  for any  volume form $dV$ on $M$ there is  a scalar function $f$ on $M$ such that
\be
{\bf PRic}_{(G, dV)} =(n-1)\Big \{f_{0|0} -f_{0}^{2}+\frac{2}{n+1}f_{0}\mathbf{S}\Big \},\label{oPRicf}
\ee
\item[(c)]  for any  volume form $dV$ on $M$ there is  a scalar function $f$ on $M$ such that
\be
{\bf Ric}_{G} =-(n-1)\Big \{\Xi_{|0} + \Xi^2 \Big \},\label{PRicf}
\ee
\een
where
 $``|"$ is the horizontal covariant derivative with respect to $\mathbf{G}$, $f_{0}=f_{x^{m}}y^{m}$,  $\Xi := \frac{{\bf S}}{n+1}-f_0 $ and ${\bf S}= {\bf S}_{(G, dV)}$.
\end{theorem}
{\it Proof}: 
(a) $\Rightarrow$ (b).   Assume that for some volume form $d\tilde{V}$,
\[ {\bf PRic}_{(G, d\tilde{V})}  = 0.\]
Then by (\ref{5297}), for any volume form  $dV= e^{-(n+1)f} d\tilde{V}$,
\begin{equation}
{\bf PRic}_{(G, dV)} = (n-1) \Big \{ f_{0|0} - f_0^2 +\frac{2}{n+1} f_0 {\bf S} \Big \} , \label{PRic*}
\end{equation}

(b) $\Rightarrow$ (c).  It follows from (\ref{eqP1}) and (\ref{PRic*}) that
\begin{eqnarray*}
{\bf Ric}_G  & = & (n-1) \Big \{ f_{0|0} - f_0^2 +\frac{2}{n+1} f_0 {\bf S} \Big \} - (n-1) \Big \{  \Big [ \frac{\bf S}{n+1}\Big ]^2 + \frac{{\bf S}_{|0}}{n+1} \Big \}\\
& = & -(n-1) \Big \{  \Xi_{|0} +\Xi^2 \Big \},
\end{eqnarray*}
where $\Xi =  \frac{\bf S}{n+1}-f_0$.

(c) $\Rightarrow$ (a).  Let $dV$ be an arbitrary  volume form on $M$.  There is a scalar function $f $ on $M$ such that (\ref{PRicf}) holds.
\[ {\bf Ric}= -(n-1) \Big \{ \Xi_{|0}+ \Xi^2 \Big \},\]
where $\Xi= \frac{\bf S}{n+1}-f_0 $ and ${\bf S}= {\bf S}_{(G, dV)}$.
Let
$d\tilde{V}= e^{(n+1) f} dV$,  we have
\[ \tilde{\bf S} = {\bf S} -(n+1) f_0 = (n+1) \Xi.\]
By (\ref{eqP2}), we get
\[ {\bf PRic} _{(G, d\tilde{V})} = {\bf Ric}_G + (n-1) \Big \{  \Xi^2 + \Xi_{|0} \Big \} =0.\]
\qed

\bigskip

As we have mentioned in the introduction, if a spray ${\bf G}$ on an $n$-dimensional manifold $M$ satisfies
that ${\bf Ric}_G=0$ and ${\bf S}_{(G, dV)} = (n+1) \phi_0 $ for some volume form $dV$ and scalar function $\phi$ on $M$, then
(\ref{PRicf}) holds. Thus ${\bf G}$ is projectively Ricci-flat.

\bigskip
\noindent{\it Proof of Theorem \ref{thmrigidit}}: Let $dV$ be an arbitrary volume form on $M$. By Theorem \ref{5310},   there is a scalar function $f$ on $M$ such that
\be
{\bf Ric}_G =- (n-1) \Big \{ \Xi_{|0} + \Xi^2 \Big \},  \label{RicXiXi}
\ee
where $\Xi :=  \frac{\bf S}{n+1}-f_0$ and ${\bf S}= {\bf S}_{(G, dV)}$.
Assume that $\Xi =\Xi (x, y) \not=0$ for some non-zero $y\in T_xM$. Let $c(t)$ be the geodesic with $c(0)=x$ and $c'(0)=y$.  By assumption on the completeness, we may assume that $c$ is defined on  $R=(-\infty, \infty)$.
Let $\Xi(t):= \Xi(c(t), c'(t))$.  Then $\Xi'(t) = \Xi_{|0} (c(t), c'(t))$. 
It follows from (\ref{RicXiXi}) that
\[  \Xi'(t) +\Xi(t)^2 \leq 0.\]

\noindent{\it Case 1}: $\Xi(0) <0$. Let $ r := \sup \{ b >0\  | \   \Xi(t) <0,  \ 0\leq t < b\}\leq +\infty$.
Observe that for $0 \leq t < r$,
\[  \Big [ \frac{1}{\Xi(t)} \Big ] ' \geq 1.\]
\[ \frac{1}{\Xi(t)} - \frac{1}{\Xi(0) }  \geq t .\]
\[ 0 >  \frac{1}{\Xi(t)} \geq  \frac{1+\Xi(0) t}{\Xi(0)}.\]
This implies that $1+\Xi(0) t >0$ for $0 \leq t <r$.  Thus $ r \leq - 1/\Xi(0)$ is finite.  Then $\Xi(r)=0$  by  the definition of $r$. Observe that
\[ \Xi (t) \leq   \frac{\Xi(0)}{1+\Xi(0) t } \leq \Xi(0) <0, \ \ \ \ 0\leq t < r.\]
Thus $\Xi(r) \leq  \Xi(0) <0$. This is impossible.

\noindent{\it Case 2}: $\Xi (0) > 0$.  Let $ r = \sup \{ b >0 \ | \ \Xi(t) > 0, \  -b < t\leq 0 \}$. 
By a similar argument for  $\Xi(t)$ on  $(-r , 0]$, we can show that this is impossible. 

Therefore $\Xi =0$. Then it follows from (\ref{RicXiXi}) that ${\bf Ric}_G =0$.
\qed

\bigskip

Below is a specific  non-trivial example.

\begin{example}
Let $\alpha_{1}=\sqrt{a_{ij}y^{i}y^{j}}$ and $\alpha_{2}=\sqrt{\bar{a}_{ij}y^{i}y^{j}}$ be two Ricci-flat Riemannian
metrics on the manifolds $M_{1}$ and $M_{2}$, respectively. Consider the following 4-th root metric
\begin{equation*}
F:=\sqrt[4]{\alpha_{1}^{4}+2c\alpha_{1}^{2}\alpha_{2}^{2}+\alpha_{2}^{4}},
\end{equation*}
where $0 < c \leq 1$ is a constant.  It is a Riemannian metric when $c =1$.
This is a Ricci-flat $(\mathbf{Ric}_F=0)$ and Berwald metric on $M:=M_{1}\times M_{2}$.  Thus for the Busemann-Hausdorff volume form $dV=dV_F$, the S-curvature ${\bf S}_{(F, dV)}=0$.
 Therefore  $F$ is projectively Ricci-flat.
\end{example}

\section{Randers metrics}
 In this section, we will
derive the equivalent conditions for a Randers metric satisfying equation\eqref{PRicf}.

Let $F=\alpha+\beta$ be a Randers metric on an $n$-dimensional manifold $M$, where $\alpha=\sqrt{a_{ij}(x)y^{i}y^{j}}$ and $\beta=b_{i}y^{i}$.  Put
\begin{eqnarray*}
r_{ij}:&=&\frac{1}{2}(b_{i;j}+b_{j;i}),~~s_{ij}:=\frac{1}{2}(b_{i;j}-b_{j;i}),\\
t_{ij}:&=&s_{il}s^{l}_{~j},~~
s_{i}:=b^{j}s_{ji},~~t_i:= b^jt_{ji}=s_ms^m_{\ i},
\end{eqnarray*}
where $``;"$ denotes the covariant derivative with respect to Levi-Civita connection of $\alpha$.
Through out this paper, we will use $``0"$ to denote the contraction with $y^{i}$. For example,
$s^{m}_{~~0}=s^{m}_{~j}y^{j},$ $t_{00}=t_{ij}y^{i}y^{j},$  $r_{00}=r_{ij}y^{i}y^{j}$.
The geodesic coefficients  of  $\mathbf{G} = {\bf G}_F$ are given by
\begin{equation}\label{5298}
G^{i}=\tilde{G}^{i}+Py^{i},~~\tilde{G}^{i}=\bar{G}^{i}+\alpha s^{i}_{~0},
\end{equation}
where $\bar{G}^{i}$ denotes the spray coefficients of $\alpha$ and $P=\frac{r_{00}-2\alpha s_{0}}{2F}$ (\cite{CS1}). Clearly, $\mathbf{G} = y^i\pxi -2 G^i\pyi$ and $\tilde{\mathbf{G}}=y^i\pxi-2\tilde{G}^i\pyi $ in \eqref{5298} are projectively equivalent. Recall that the projective spray $\hat{\bf G} = {\bf G} +\frac{2{\bf S}}{n+1} {\bf Y}$  is   projectively invariant for a fixed volume form $dV=\sigma (x) dx^1 \cdots dx^n $. Thus the projective Ricci curvature of  $({\mathbf{G}},dV)$ is given by
\begin{equation}\label{5292}
{\bf PRic}_{(G, dV)}  = \mathbf{PRic}_{(\tilde{\mathbf{G}},dV)}=\mathbf{Ric}_{\tilde{\mathbf{G}}}+ (n-1) \Big \{ \frac{\tilde{\mathbf{S}}_{|0}}{n+1}
+\Big [ \frac{\tilde{\mathbf{S}}}{n+1} \Big ]^{2}\Big \},
\end{equation}
where $\tilde{\mathbf{S}}={\bf S}_{(\tilde{G}, dV)} $ is the $S$-curvature of  $(\tilde{\mathbf{G}},dV)$ and $``|"$ denotes horizontal covariant derivative with respect to $\tilde{\mathbf{G}}$.

\bigskip
\noindent
{\it Proof of Theorem \ref{6183}}:
We shall calculate each term on the right side of (\ref{5292}) as follows. Firstly, we have known that
\begin{equation}\label{5293}
\mathbf{Ric}_{\tilde{\mathbf{G}}}={\rm Ric}_{\alpha}+(2\alpha s^{m}_{~0;m}-2t_{00}-\alpha^{2}t^{m}_{~m}),
\end{equation}
where $``;"$ denotes the covariant derivative with respect to Levi-Civita connection of $\alpha$ and ${\rm Ric}_{\alpha}$ is the Ricci curvature of $\alpha$.
From $\tilde{G}^{i}=\bar{G}^{i}+\alpha s^{i}_{~0}$, it follows that
\begin{eqnarray}
\nonumber\tilde{\mathbf{S}}&=&\frac{\partial \tilde{G}^{m}}{\partial y^{m}}-y^{m}\frac{\partial}{\partial x^{m}}(\ln \sigma)\\
\nonumber&=&\frac{\partial \bar{G}^{m}}{\partial y^{m}}-y^{m}\frac{\partial}{\partial x^{m}}(\ln \sigma_{\alpha} )+ y^{m}\frac{\partial}{\partial x^{m}}\Big (\ln \frac{\sigma_{\alpha}}{\sigma}\Big )\\
\label{5294}&=&y^{m}\frac{\partial}{\partial x^{m}}\Big (\ln \frac{\sigma_{\alpha}}{\sigma}\Big ).
\end{eqnarray}
Here we have used $\frac{\partial \bar{G}^{m}}{\partial y^{m}}-y^{m}\frac{\partial}{\partial x^{m}}\ln \sigma_{\alpha}=0$ and
$\frac{\partial}{\partial y^{m}}(\alpha s^{m}_{~0})=0$, where  $\sigma_{\alpha}=\sqrt{\text{det}(a_{ij}(x))}$. Denote by $\mu=\frac{1}{n+1}\ln \frac{\sigma_{\alpha}}{\sigma}$ and $\mu_{0}=\mu_{x^{m}}y^{m}=\frac{\partial \mu}{\partial x^{m}}y^{m}$. Then (\ref{5294}) can be rewritten by $\tilde{\mathbf{S}}=(n+1)\mu_{0}$. Thus we have
\begin{eqnarray}
\nonumber\tilde{\mathbf{S}}_{|0}&=&\tilde{\mathbf{S}}_{;0}-2\alpha s^{m}_{~0}\tilde{\mathbf{S}}_{\cdot m}\\
\label{5295}&=&(n+1)(\mu_{0;0}-2\alpha s^{m}_{~0}\mu_{x^{m}}).
\end{eqnarray}
and
\begin{eqnarray}\label{5296}
\tilde{\mathbf{S}}^{2}&=&(n+1)^{2}\mu_{0}^{2}.
\end{eqnarray}
Plugging (\ref{5293}), (\ref{5295}) and (\ref{5296}) into (\ref{5292}) yields
\begin{eqnarray}
\nonumber\mathbf{PRic}_{(G, dV)} &=&{\rm Ric}_{\alpha}+(2\alpha s^{m}_{~0;m}-2t_{00}-\alpha^{2}t^{m}_{~m})\\
\label{5299}&+&(n-1)(\mu_{0;0}-2\alpha s^{m}_{~0}\mu_{x^{m}})+(n-1)\mu_{0}^{2}.
\end{eqnarray}
On the other hand, by (\ref{5298}), the $S$-curvature of $(F,dV)$ is given by
\begin{eqnarray}
\nonumber\mathbf{S}&=&\frac{\partial G^{m}}{\partial y^{m}}-y^{m}\frac{\partial}{\partial x^{m}}(\ln \sigma)\\
\nonumber&=&\frac{\partial \bar{G}^{m}}{\partial y^{m}}+(n+1)P-y^{m}\frac{\partial}{\partial x^{m}}(\ln \sigma)\\
\nonumber&=&y^{m}\frac{\partial}{\partial x^{m}}(\ln \sigma_{\alpha})+(n+1)P-y^{m}\frac{\partial}{\partial x^{m}}(\ln \sigma)\\
\label{5300}&=&(n+1)(\mu_{0}+P),
\end{eqnarray}
where $P=\frac{r_{00}-2\alpha s_{0}}{2F}$. Hence by (\ref{5299}) and (\ref{5300}), equation (\ref{oPRicf}) is equivalent to
\begin{eqnarray*}
&&{\rm Ric}_{\alpha}+(2\alpha s^{m}_{~0;m}-2t_{00}-\alpha^{2}t^{m}_{~m})\\
&+&(n-1)\{\mu_{0;0}-2\alpha s^{m}_{~0}\mu_{x^{m}}\}+(n-1)\mu_{0}^{2}\\
&=&(n-1)\{-f_{0}^{2}+f_{0;0}-2Pf_{0}-2\alpha s^{m}_{~0}f_{x^{m}}+\frac{2}{n+1}f_{0}[(n+1)(\mu_{0}+P)]\}\\
&=&(n-1)\{-f_{0}^{2}+f_{0;0}-2\alpha s^{m}_{~0}f_{x^{m}}+2f_{0}\mu_{0}\}.
\end{eqnarray*}
The above equation is actually in the form of
\begin{equation*}
A+\alpha B=0,
\end{equation*}
where $A$ and $B$ are polynomial in $y$. We see that this equation is valid if and only if $A=B=0$.
Then we obtain
\begin{eqnarray*}
{\rm Ric}_{\alpha} & = & 2 t_{00}+\alpha^2 t^m_m 
-(n-1) \{ \mu_{0;0}  - f_{0;0}+(\mu_0-f_0)^2 \},\\
s^m_{\ 0;m} & = & (n-1) s^m_{\ 0}\{  \mu_{x^m} -  f_{x^m}\}.
\end{eqnarray*}
Letting $h:= \mu - f$,  we complete the proof of Theorem \ref{6183}.
\qed

\section{The Riemannian case}
Let $F=\alpha$ be a Riemann metric  and $dV=e^{-(n+1) f}dV_{\alpha}$ be a volume form on $M$.  Then
${\bf S}_{(G, dV)} = (n+1) f_{;0}$.
By (\ref{630}), we have
\[ {\bf PRic}_{(G, dV)} = {\bf Ric}_{\alpha}  + (n-1) \Big \{ f_{0;0} + f_{;0}^2 \Big \}
.\]
Then Theorem \ref{6183} reduces to the following
\begin{cor}\label{6152}
Let $(M,\alpha)$ be an $n$-dimensional Riemannian manifold and $dV=e^{-(n+1) f}dV_{\alpha}$ a volume form on $M$. Then
${\bf PRic}_{(G, dV)} = 0$
if and only if
\begin{equation*}
\mathbf{Ric}_{\alpha}=-(n-1)\{f_{0;0}+f_{;0}^{2}\},
\end{equation*}
\end{cor}

This leads to a new notion of {\it weighted Ricci curvature} ${\bf Ric}^f_{\alpha} $ on a Riemannian  $n$-manifold $(M, \alpha)$ with a scalar function $f=f(x)$ on $M$:
\begin{equation}
  {\bf Ric}^f_{\alpha} := {\bf Ric}_{\alpha} + (n-1) \{ f_{0;0} +f_{;0}^2 \}.
\end{equation}
This weighted Ricci curvature and its relationship with other geometric quantities of $(\alpha, f)$ deserves further study.

\begin{example}
Let $(\Bbb{S}^{n},\alpha)$ be a Riemannian $n$-sphere with constant sectional curvature $\mathbf{K}_{\alpha}=1$.
Recall that in this case, there exists a scalar function $\phi$ on $M$ such that $\phi_{0;0}=-\phi\alpha^{2}$, where $``;"$ is covariant derivative
with respect to the Riemann metric $\alpha$. Thus, the Ricci curvature ${\mathbf{Ric}}_{\alpha}$ of $\alpha$ satisfies
\begin{equation}\label{6151}
{\mathbf{Ric}}_{\alpha}=(n-1)\alpha^{2}=-(n-1)\frac{\phi_{0;0}}{\phi}.
\end{equation}
Set $f=\ln |\phi|$ which is only defined on $\{ \phi \not=0\}$. Then (\ref{6151}) can be rewritten by
\begin{equation*}
{\mathbf{Ric}}_{\alpha}+(n-1)\{f_{0}^{2}+f_{0;0}\} =0.
\end{equation*}
It is the case in Corollary \ref{6152}.
\end{example}

\vspace{0.6cm}

\noindent Zhongmin Shen

\noindent Department of Mathematical Sciences,

\noindent  Indiana University-Purdue University Indianapolis, IN 46202-3216, USA.

\noindent \verb"zshen@iupui.edu"

\vspace{0.6cm}

\noindent Liling Sun

\noindent Department of Mathematics,

\noindent  Taiyuan University of Technology Taiyuan, 030024, P.R.China.

\noindent \verb"sunliling@yeah.net"

\end{document}